\documentclass{article}

\usepackage[english,russian]{babel}
\usepackage[tbtags]{amsmath}
\usepackage{amsfonts,amssymb}
\usepackage[usenames]{color}

\usepackage{graphicx}

\usepackage[matrix,arrow,curve]{xy}
\usepackage{lscape}

\numberwithin{equation}{subsection}

\newtheorem{theorem}{Теорема}
\newtheorem{definition}{Определение}

\newtheorem{lemma}{Лемма}

\newtheorem{sled}{Следствие}
\newtheorem{suggestion}{Предложение}

\newtheorem{proof}{Доказательство}
\begin{document}
\righthyphenmin=2

\author{А.\,М.~Михович}
\title{Квазирациональность и асферические (про-$p$)копредставления}
\markboth{А.\,М.~Михович}{Квазирациональность и асферические копредставления}
\maketitle
\begin{abstract}
   В работе показано, что свойство квазирациональности (про-$p$) копредставления (про-$p$) группы $G$ является свойством самой (про-$p$) группы и не зависит от выбора копредставления. Доказано, что класс квазирациональных копредставлений шире класса асферических про-$p$-копредставлений (комбинаторно-асферических копредставлений в дискретном случае). Для квазирациональных копредставлений введено понятие обобщенной пермутационности модуля соотношений, которое оказывается равносильным пермутационности его $mod(p)$ фактора.

  Библиография: 22 названия.

  Ключевые слова: квазирациональность, асферичность, обобщенная пермутационность.
\end{abstract}
\footnotetext{}

\subsection{Введение}
В работе \cite{Mikh2014} мы ввели понятие квазирационального (про-$p$) копредставления конечного типа, которое в дискретном случае включает в себя асферические копредставления и их подкопредставления, а также про-$p$-группы с одним соотношением в про-$p$-случае.
Квазирациональные копредставления естественно обобщают дискретные комбинаторно-асферические копредставления (введенные Й.Хьюбшманном \cite{Hue,Hue1,CCH}), асферические (про-$p$) копредставления (изученные О.В.Мельниковым в \cite{Mel1}) и удобны тем, что позволяют единообразно описывать модули соотношений про-$p$-групп с помощью техники аффинных групповых схем \cite{Mikh2014,Mikh2015,Mikh2016}. В свете знаменитой гипотезы ``асферичности'' Уайтхеда крайне любопытно выяснить различие между квазирациональностью и асферичностью (удобнее сравнивать с комбинаторной асферичностью), которое действительно имеет место и проявляется в доказательстве Предложения \ref{p002} на копредставлениях конечных $p$-групп (в частности показано, что такими являются копредставления группы кватернионов $Q_8$).

В разделе 2 предлагаемой работы показано, что введенное ранее определение \cite[Определение 1]{Mikh2014} квазирационального (про-$p$)копредставления (про-$p$)группы $G$ эквивалентно утверждению о том, что группы гомологий $H_2(G,\mathbb{Z})$ (соответственно $H_2(G,\mathbb{Z}_p)$ в про-$p$-случае) не имеют кручения. Таким образом свойство квазирациональности (про-$p$) копредставления (про-$p$)группы $G$ является свойством самой (про-$p$)группы и (в отличии от свойства копредставления дискретной группы быть асферическим) не зависит от выбора копредставления. Итак, группа $Q_8$ это пример квазирациональной, но не асферической про-$p$-группы (в \cite{Mel1} доказано, что асферическая конечная $p$-группа циклична).



В разделе \ref{s3} с помощью нового понятия \emph{обобщенной пермутационности} мы ``унифицируем'' введенные в \cite{Mel1} понятия $\mathbb{Z}_p$ и $\mathbb{F}_p$- пермутационности, которые являются квинтэссенциями понятий $CA$-асферичности в дискретном случае (поскольку тут данные понятия обобщают $CA$-асферичность) и асферичности в смысле Мельникова для про-$p$-копредставлений (в \cite{Mel1} показано, что понятие $\mathbb{F}_p$- пермутационности эквивалентно про-$p$-аналогу $CA$-асферичности). В Теореме \ref{t0} доказана эквивалентность обобщенной пермутационности модуля соотношений и пермутационности его $mod(p)$ фактора для $QR$- копредставлений (как поясняется в разделе \ref{s3} надеяться на совпадение данных понятий в общем случае было бы слишком оптимистично). Показано, что класс $QR$-копредставлений соответствует ожиданиям О.В. Мельникова, которые мы формулируем в виде Следствия \ref{H1} (``о существовании оболочки''). Отметим, что О.В.Мельников и его ученики надеялись выделить подходящий класс про-$p$-групп с помощью когомологических данных (планарные про-$p$-группы \cite{MS}, теория концов групп и виртуальная двойственность Пуанкаре \cite{Kor}), но реализовать желанные свойства класса из Следствия \ref{H1} в полном объеме не получилось.

Автор выражает глубокую благодарность А.С. Мищенко и В.М. Мануйлову за их внимание, критику, бесконечное терпение и помощь автору данной работы.

\subsection{Квазирациональные копредставления}
\label{subsec1}
Копредставление конечного типа дискретной группы $G$ - это точная последовательность

\begin{equation}
1 \rightarrow R \rightarrow F \xrightarrow{\pi} G \rightarrow 1 \label{eq1}
\end{equation}
в которой $F=F(X)$ есть свободная группа с конечным множеством образующих $X$, а  $R$ - нормальная подгруппа в $F$, порожденная конечным числом определяющих соотношений $r\in R$.
Про-$p$-группой называется группа, изоморфная проективному пределу конечных $p$-групп. Это топологическая группа (с топологией проективного предела), которая является компактной вполне несвязной группой. Для таких групп имеется теория копредставлений во многом аналогичная комбинаторной теории дискретных групп \cite{Koch}, \cite{ZR}.
Пусть $X$ - проконечное топологическое пространство, тогда в категории про-$p$-групп имеется свободная про-$p$-группа $(F(X), i:X\hookrightarrow F(X)),$ обладающая свойствами свободного объекта. А именно, для любого непрерывного гомоморфизма $\phi:X\rightarrow G$ в про-$p$-группу $G$ такого, что образ $X$ порождает $G,$ существует и единственный непрерывный гомоморфизм про-$p$-групп $\widetilde{\phi}:F(X)\rightarrow G,$ делающий соответствующую диаграмму коммутативной \cite[3.3]{ZR}. Если $X$ - конечное множество, то свободная про-$p$-группа может быть построена конструктивно, как $F(X)=\varprojlim_{U\triangleleft \Phi(X)}\Phi(X)/U,$ где $\mid\Phi(X)/U\mid=p^n$ - про-$p$-пополнение дискретной свободной группы $\Phi(X)$ (порожденной образующими из $X$) в про-$p$-топологии на $\Phi(X),$ в которой базис системы окрестностей единицы составляют нормальные подгруппы конечного индекса равного некоторой степени простого числа $p$. По аналогии с копредставлением конечного типа дискретной группы, будем говорить, что про-$p$-группа $G$ задана про-$p$- копредставлением конечного типа, если группа $G$ включена в точную последовательность \eqref{eq1}, в которой $F$ - свободная про-$p$-группа с конечным числом образующих (то есть $X$ - конечное множество), а $R$ - замкнутая нормальная подгруппа, топологически нормально порожденная конечным числом элементов в $F$. При этом в теории про-$p$-групп принято рассматривать про-$p$-копредставления, у которых число образующих в $F$ совпадает с числом образующих в $G$.

Пусть $C=\varprojlim C_{\alpha}$ - проконечное кольцо ($C_{\alpha}$ - конечные кольца), тогда обозначим через $CG$ - пополненную групповую алгебру про-$p$-группы $G$. Под пополненной групповой алгеброй мы понимаем топологическую алгебру $CG=\varprojlim CG_{\mu}$ \cite[5.3]{ZR}, где $G=\varprojlim G_{\mu}$ - разложение про-$p$-группы $G$ в проективный предел конечных $p$-групп $G_{\mu}$.

Для дискретных групп, $p$ будет пробегать все простые числа, для про-$p$- групп $p$ фиксировано. Пусть $G$ - (про-$p$)группа с (про-$p$)копредставлением конечного типа \eqref{eq1}, $\overline{R}=R/[R,R]$ соответствующий \emph{$G$-модуль соотношений}, где $[R,R]$ - это коммутант, а действие $G$ индуцировано сопряжением $F$ на $R$. Для каждого простого числа $p\geq2$ обозначим через ${\Delta}_p$ - аугментационный идеал кольца $\mathbb{F}_pG,$ а через $\mathcal{M}_n, n\in \mathbb{N}$ его $p$-фильтрацию Цассенхауза в $F$ с коэффициентами в поле $\mathbb{F}_p$, которая определена по правилу
$\mathcal{M}_n=\{f \in F\mid f-1 \in {\Delta}^n_p\}.$
В про-$p$-случае под ${\Delta}^n$ мы понимаем замыкание модуля порожденного $n$-ми степенями элементов из $\Delta={\Delta}_p$, а в дискретном случае - это $n$-ая степень идеала ${\Delta}_p$ \cite{Pas}. Свойства этой фильтрации в про-$p$-случае изложены в \cite[7.4]{Koch}, в дискретном случае свойства фильтрации Цассенхауза аналогичны \cite[Гл.11]{Pas}, разница состоит в использовании обычного группового кольца вместо пополненного.

\begin{definition}\cite[Определение 1]{Mikh2014} \label{d10} Копредставление конечного типа \eqref{eq1} будем называть квазирациональным ($QR$-копредставлением), если для каждого $n>0$ и для каждого простого $p\geq2$ $F/R\mathcal{M}_n$-модуль $R/[R,R\mathcal{M}_n]$ не имеет нетривиального $p$-кручения ($p$ фиксировано для про-$p$-групп и пробегает все простые числа $p\geq2$ и соответствующие $p$-фильтрации Цассенхауза в дискретном случае). Модули соотношений таких копредставлений будем называть квазирациональными модулями соотношений.
\end{definition}

Несмотря на существование более простых эквивалентных формулировок квазирациональности (Предложение \ref{p03} и нункт 1 Предложения \ref{p002} ), именно Определение \ref{d10} открывает возможность использования аффинных групповых схем для исследования такого рода копредставлений. Дело в том, что отсутствие $p$-кручения в факторах $R/[R,R\mathcal{M}_n]$ и точность слева фуктора проективного предела обеспечивают вложение проективного предела $p$-адических пополнений $\varprojlim (R/[R,R\mathcal{M}_n])^{\wedge}_p$ в рационализированный модуль соотношений $\overline{R} \widehat{\otimes}\mathbb{Q}_p=\varprojlim_n R/[R,R \mathcal{M}_n]\otimes \mathbb{Q}_p,$ где $\mathbb{Q}_p$ - поле $p$-адических чисел.
В про-$p$-случае последний отождествляется с $\overline{R^{\wedge}_w} (\mathbb{Q}_p)$ \cite[Лемма 2]{Mikh2016}, это так называемая абелианизация $p$-адического непрерывного проунипотентного пополнения $R$. Последнее наблюдение позволяет включить $\overline{R} \widehat{\otimes}\mathbb{Q}_p$ в коммутативную диаграмму \cite[Теорема 1 и Следствие 2]{Mikh2016} и использовать методы теории когомологий проунипотентных групп. В частности, в работе \cite[Теорема 1]{Mikh2017} представлен аналог ``теоремы о тождествах'' для про-$p$-групп c одним соотношение, а в \cite[Предложение 3]{Mikh2016} удается получить достаточное условие когомологической размерности 2 для про-$p$-групп с одним соотношением.

\begin{suggestion}\cite[Предложение 4]{Mikh2014} \label{p03} Пусть задано (про-$p$) копредставление конечного типа \eqref{eq1}, тогда следующие свойства равносильны:

(i) \eqref{eq1} является $QR$-(про-$p$)копредставлением;

(ii) (про-$p$)фактор-модуль коинвариантов $R/[R,F]$ не имеет кручения.
\end{suggestion}

Отметим, что квазирациональность можно было бы определить с помощью любого базиса системы окрестностей единицы свободной (про-$p$)группы, состоящего из нормальных делителей индекса равного степени заданного простого числа $p,$ однако фильтрация Цассенхауза возникает из алгебры Хопфа группового кольца над полем $\mathbb{F}_p,$ что гармонично вписывается в концепцию схематизации, используемую в \cite{Mikh2014,Mikh2015,Mikh2016}.

Пусть \eqref{eq1} - копредставление дискретной группы, обозначим $K(X;R)$ его стандартный двумерный $CW$-комплекс (детали в \cite{CCH}). Мы называем копредставление \eqref{eq1} \textbf{асферическим}, если $K(X;R)$ является асферическим, что эквивалентно условию $\pi_q(K(X;R))=0$, $q\geq 2$.

Отметим, что представленные в данной работе алгебраические конструкции имеют геометрическую природу. Действительно, выберем в качестве универсального накрытия $K(X;R)$ комплекс Кэли $\widetilde{K(X;R)}$ и рассмотрим \cite[I,Теорема 5.3]{Bro} цепной комплекс свободных $\mathbb{Z}G$-модулей его одномерного остова $Y$: $$C_1Y\rightarrow C_0Y\rightarrow \mathbb{Z}\rightarrow 0.$$ Тогда $H_1 Y\cong (\pi_1 Y)_{ab}\cong \overline{R}$ и мы имеем дело с копредставлениями с особой $G$-модульной структурой на $H_1 Y.$

Пусть $r\in F$ - некоторый элемент дискретной свободной группы $F$, обозначим через $\sqrt{r}$ корень элемента в $F$, по определению - это некоторый элемент $s\in F,$ что $s^m=r$ и $m$ максимально (так как $F$ - свободная группа, то такой элемент существует и единственен). Обычная асферичность является в определенном смысле слишком узким понятием и поэтому нуждается в расширении. Имеется привлекательное понятие комбинаторной асферичности (\textbf{CA-асферичности}), которое может быть определено следующим образом

\begin{definition}\cite[1]{CCH} \label{d001}
Копредставление \eqref{eq1} группы $G$ называется $CA$- асферическим, если модуль соотношений $\overline{R}$ раскладывается как $\mathbb{Z}G$-модуль в прямую сумму циклических (то есть однопорожденных) подмодулей $P_r,r\in R$, где $P_r$ порождается $\overline{r}=r[R,R]$ с единственным соотношением $\pi(s)\overline{r}=\overline{r},s=\sqrt{r}$, а $\pi$ - гомоморфизм из \eqref{eq1}.
\end{definition}

Обычные асферические копредставления это те $CA$- асферические копредставления, у которых все циклические подмодули $P_r,r\in R$ в разложении модуля соотношений из Определения \ref{d001} изоморфны свободному модулю $\mathbb{Z}G$, то есть корни из определяющих соотношений совпадают с самими определяющими соотношениями \cite{BH}. Иначе говоря, $CA$-асферичность означает, что у данного копредставления нет нетривиальных "сизигий" (тождеств между соотношениями).

Важность класса $QR$-(про-$p$)копредставлений обусловлена тем, что он имеет непосредственное отношение к ряду открытых задач двумерной топологии, таких как знаменитая гипотеза $"$асферичности$"$ Уайтхеда \cite{BH}, а также к проблемам комбинаторной теории про-$p$-групп, например к вопросу Серра о структуре модуля соотношений про-$p$-групп с одним соотношением \cite[10.2]{Se}.
Напомним, что $QR$--копредставления включают в себя $CA$-асферические копредставления и все их подкопредставления в дискретном случае \cite[Предложение 3]{Mikh2014} ($R/[R,F]\cong \oplus_r (P_r)_F\cong \oplus_r \mathbb{Z}$ и коинварианты по действию $F$ не содержат кручения), а также про-$p$-копредставления про-$p$-групп с одним соотношением \cite[Предложение 1]{Mikh2014}. Введем обозначение $\mathbb{Z}_{(p)}$ для $\mathbb{Z}$ в случае дискретных групп и для $\mathbb{Z}_{p}$ в случае про-$p$-групп.

\begin{suggestion} \label{p002}
1) (Про-$p$)копредставление \eqref{eq1} (про-$p$)группы $G$ является квазирациональным тогда и только тогда, когда $H_2(G,\mathbb{Z}_{(p) })$ не имеет кручения.

2) Класс дискретных $QR$-копредставлений строго шире класса $CA$- асферических копредставлений.

\end{suggestion}
\begin{proof}
1) Рассмотрим длинную точную последовательность групп гомологий для копредставления \eqref{eq1}
$$0\rightarrow H_1(G,IG)\rightarrow R/[R,F]\rightarrow \mathbb{Z}_{(p)}^{|X|}\rightarrow (IG)_G\rightarrow 0,$$ которая получается из короткой точной последовательности Кроуэлла-Линдона \cite[Теорема 2.2]{I} $\mathbb{Z}_{(p) }G$- модулей
$$0\rightarrow \overline{R}\rightarrow \mathbb{Z}_{(p) }G^{|X|}\rightarrow IG\rightarrow 0,$$ где $IG$ - аугментационный идеал $\mathbb{Z}_{(p)}G$, а через $(IG)_G$ мы как обычно обозначили фактор-модуль $G$-коинвариантов $G$-модуля $IG$. Теперь воспользуемся изоморфизмом $H_1(G,IG)\cong H_2(G,\mathbb{Z}_{(p) }),$ который следует из сдвига размерности, возникающего из точной последовательности аугментации $$0\rightarrow IG\rightarrow \mathbb{Z}_{(p)}G\rightarrow \mathbb{Z}_{(p)}\rightarrow 0.$$ Откуда видно, что подгруппа кручения в $R/[R,F]$ изоморфна подгруппе кручения в  $H_2(G,\mathbb{Z}_{(p) })$. Из Предложения \ref{p03} следует, что квазирациональность эквивалентна отсутствию кручения в фактор-группе коинвариантов $R/[R,F]$, откуда получаем первое утверждение.

2) Второе утверждение следует из существования конечных нециклических $p$- групп с тривиальным мультипликатором Шура \cite[3.4]{Kar} (такие конечные $p$-группы бывают только \emph{с двумя и тремя образующими}, например группы кватернионов $Q_n, n\geq1$ \cite[2.4.8]{Kar}). На самом деле, для конечных групп мультипликатор Шура изоморфен $H_2(G,\mathbb{Z})$ \cite[2.7.3]{Kar}, а поэтому всякая $p$-группа c тривиальным мультипликатором Шура квазирациональна. Ввиду результатов \cite{Hue}, когомологии $H^n(G,\mathbb{Z})$ CA- асферических групп периодичны при $n\geq3$ с периодом 2 \cite[Теорема 2]{Hue} (поскольку таков период у когомологий циклических групп), но например у $Q_8$ когомологии периодичны с периодом 4 \cite[41.1]{Jo}, поэтому у $$Q_8=\langle a,b| a^4=b^2,aba=b\rangle$$нет CA-асферических копредставлений и следовательно класс $QR$- копредставлений строго шире.
\end{proof}

Если про-$p$-копредставление конечной $p$-группы $G$ является асферическим (в случае дискретного копредставления следует потребовать минимальность), то, как показывает Мельников \cite[Теорема 2.7]{Mel1}, такая конечная $p$-группа является циклической $G\cong \mathbb{Z}/p^n\mathbb{Z}$ (в дискретном случае рассуждения Мельникова тоже остаются в силе, но избавиться от требования минимальности копредставления мы не можем).
Собирая предыдущие рассуждения и исправляя неточность в \cite[Замечание 1]{Mikh2014}, ставшую катализатором данной работы, получаем, что всякое (минимальное) копредставление конечной нециклической $p$-группы $G$ c тривиальным мультипликатором Шура является квазирациональным, но не ($CA$)асферическим.

\begin{sled}
Cвойство квазирациональности не зависит от выбора копредставления (про-$p$) группы $G$ и является свойством самой (про-$p$) группы.
Таким образом, если одно из копредставлений $G$ квазирационально, то и все другие тоже являются квазирациональными.
\end{sled}

\subsection{Обобщенная пермутационность} \label{s3}

Свободным $\mathbb{F}_p$-модулем, где $\mathbb{F}_p$ - поле из $p$ элементов, над пунктированным проконечным пространством $(T,t_0)$ \cite[1.7]{Mel1} называется проконечный $\mathbb{F}_p$-модуль $\mathbb{F}_p(T,t_0)$ такой, что существует вложение $\omega:T\rightarrow \mathbb{F}_p(T,t_0),\omega(t_0)=0,$  обладающее следующим универсальным свойством:

(*) для любого непрерывного вложения $\gamma:T\rightarrow B$ такого, что $\gamma(t_0)=0,$ где $B$- проконечный $\mathbb{F}_p$-модуль, существует единственный гомоморфизм проконечных $\mathbb{F}_p$-модулей $\alpha:\mathbb{F}_p(T,t_0)\rightarrow B,$ удовлетворяющий равенству $\gamma=\alpha\omega.$

В этой части мы работаем со cвободными $\mathbb{F}_p$-модулями $\mathbb{F}_p(S)$ над не пунктированными проконечными пространствами $S,$ которые возникают как проективные пределы конечномерных $\mathbb{F}_p(S_{\alpha})$, но данная конструкция является частным случаем предыдущего определения, поскольку достаточно взять пунктированное пространство $T=S\cup\{t_0\}б,$ где точка $t_0$ изолирована в пространстве $T$.

Пусть теперь $G$ - про-$p$-группа, тогда будем говорить, что $(T,t_0)$ - это $G$ - пространство, если $G$ действует непрерывно на проконечном пространстве $(T,t_0)$ гомеоморфизмами, оставляющими на месте отмеченную точку $t_0$.

\begin{definition} \cite[1.8]{Mel1}
Модуль $\mathbb{F}_p(T,t_0),$ где $(T,t_0)$ какое-либо $G$- пространство, называется пермутационным $G$-модулем, если действие каждого элемента $g\in G$ является автоморфизмом этого модуля, продолжающим по универсальному свойству (*) непрерывное отображение $t\mapsto g\cdot t$ из $T$ в $T\subset\mathbb{F}_p(T,t_0).$
\end{definition}

Копредставления дискретных группы с одним соотношением ввиду теоремы Линдона $"$о тождествах$"$ - базовый пример $CA$-асферческих копредставлений. Параллельно работам \cite{CCH,Hue,Hue1} О.В.Мельников ввел в рассмотрение и изучил \cite{Mel1} (с прицелом на про-$p$- группы с одним соотношением) копредставления про-$p$-групп, у которых $\mathbb{F}_p$- модуль соотношений является пермутационным.

\begin{definition} $\mathbb{F}_pG$-модуль соотношений $\overline{R}/p\overline{R}$ копредставления \eqref{eq1} называется $\mathbb{F}_p$- пермутационным, если имеется изоморфизм пунктированных проконечных $G$-модулей $\overline{R}/p\overline{R}=R/R^p[R,R]\cong \mathbb{F}_p(T,t_0),$ где $(T,t_0)$- некоторое проконечное $G$-пространство c отмеченной точкой.
\end{definition}

В полной аналогии с понятием $\mathbb{F}_p$-пермутационного модуля соотношений $\overline{R}/p\overline{R}$ введем понятие $\mathbb{Z}_p$-пермутационного модуля соотношений $\overline{R}$.

\begin{definition} $\mathbb{Z}_pG$-модуль соотношений $\overline{R}$ копредставления \eqref{eq1} называется $\mathbb{Z}_p$-пермутационным, если имеется изоморфизм пунктированных проконечных $G$- модулей $\overline{R}=R/[R,R]\cong \mathbb{Z}_p(T,t_0),$ где $(T,t_0)$- некоторое проконечное $G$-пространство c отмеченной точкой.
\end{definition}

Пусть $G$ - дискретная группа, обозначим через $G^{\wedge}_p$ ее про-$p$ пополнение, которое по определению задается, как $G^{\wedge}_p=\varprojlim_{N\in \mathcal{N}} G/N,$ где $\mathcal{N}$ это система нормальных делителей конечного индекса равного степени $p$, то есть $\mathcal{N}=\{N\unlhd G, \mid G/N\mid=p^n, n\in \mathbb{N}\} \cite[2.1.6]{ZR}.$ Если $G$ про-$p$-группа, то считаем, что $G^{\wedge}_p=G.$ Каждому (про-$p$) модулю соотношений, возникающему из \eqref{eq1}, поставим в соответствие две проективные системы:
$$\mathcal{T}_n=((\frac{R} {[R,R\mathcal{M}_n]} )^{\wedge}_p, \phi^{n+1}_{n}),\mathcal{T}_n^p=(\frac{R} {R^p[R,R\mathcal{M}_n]}, \widetilde{\phi}^{n+1}_{n}),$$
где $n\in \mathbb{N},$ с $G^{\wedge}_p$-модульными гомоморфизмами  $$\phi^{n+1}_{n}:(R/[R,R\mathcal{M}_{n+1}])^{\wedge}_p\twoheadrightarrow (R/[R,R\mathcal{M}_{n}])^{\wedge}_p,$$ $$\widetilde{\phi}^{n+1}_{n}:R/R^p[R,R\mathcal{M}_{n+1}]\twoheadrightarrow R/R^p[R,R\mathcal{M}_{n}],$$ индуцированными вложениями $[R,R\mathcal{M}_{n+1}]\subseteq [R,R\mathcal{M}_{n}]$.

Пусть $\zeta$- примитивный корень из единицы порядка, равного степени $p$ и $H$ - некоторая подгруппа в конечной $p$-группе $G$. Предположим, что $H$ действует слева на $\mathbb{Z}_p[\zeta]$ посредством гомоморфизма в группу корней из единицы $\xi:H\rightarrow \langle \zeta \rangle$, мы можем продолжить это действие по линейности до левого действия $\mathbb{Z}_p[\zeta]H$ на $\mathbb{Z}_p[\zeta]$. Разумеется, что $H$ действует справа обычным умножением на $G$, а это действие можно продолжить по линейности до действия $\mathbb{Z}_p[\zeta]H$ на $\mathbb{Z}_p[\zeta]G.$
\emph{Главным обобщенным пермутационным $G$-модулем}
$$\mathbb{Z}_p[\zeta]\uparrow^G_{H}:=\mathbb{Z}_p[\zeta]G\otimes_{\mathbb{Z}_p[\zeta]H} \mathbb{Z}_p[\zeta]$$
называется левый $\mathbb{Z}_p[\zeta]G$-модуль, индуцированный с $\mathbb{Z}_p[\zeta]H$-модуля $\mathbb{Z}_p[\zeta]$.

Обозначим $\pi=1-\zeta,$ тогда $\pi$ порождает $(\pi)$ - простой идеал в $\mathbb{Z}_p[\zeta],$ лежащий над $p\mathbb{Z}_p$ (то есть $(\pi)\cap \mathbb{Z}_p=(p)=p\mathbb{Z}_p$). Целое расширение $\mathbb{Z}_p[\zeta]$ над $\mathbb{Z}_p$ вполне разветвлено с индексом ветвления $p-1,$ а поэтому вложение $\mathbb{Z}_p\hookrightarrow\mathbb{Z}_p[\zeta]$ индуцирует изоморфизм $\mathbb{F}_p\cong\mathbb{Z}_p[\zeta]/\pi\mathbb{Z}_p[\zeta]$ \cite[Предл. 7.13]{Neu}.

\begin{definition} \label{d01} \cite[1]{Weiss} Пусть $G$ - конечная $p$-группа, тогда будем говорить, что $M$ - это обобщенный пермутационный $G$-модуль , если:

(i) $M$ - свободный конечнопорожденный $\mathbb{Z}_p[\zeta]$-модуль;

(ii) $M\cong \oplus_{i\in I}\mathbb{Z}_p[\zeta]\uparrow^G_{H_i} $, то есть $M$ изоморфен прямой сумме главных обобщенных пермутационных модулей.
\end{definition}

\begin{definition} \label{d1} Пусть $G$ - про-$p$-группа и фиксирован некоторый базис системы окрестностей единицы $\mathfrak{U}$ в $G$, состоящий из открытых нормальных делителей конечного индекса $U\lhd G$. Предположим, что задана $M=\{ M_U\}_{U\in\mathfrak{U}}$ - некоторая проективная система, состоящая из $\mathbb{Z}_p$- проективных $\mathbb{Z}_p[G/U]$-модулей конечной $\mathbb{Z}_p$- размерности. Будем говорить, что $M$ - это обобщенный пермутационный про-$G$-модуль, если для каждой $U\in \mathfrak{U}$ существует гомоморфизм $\xi_U:G/U\rightarrow \langle\zeta\rangle,$ где $\zeta$ - примитивный корень $p$-й степени из единицы в $\mathbb{Z}_p$, что $G/U$-модуль $M_U$ является обобщенным пермутационным $G/U$-модулем в смысле Определения \ref{d01}.
\end{definition}

  В случае про-$p$-копредставлений \eqref{eq1} $M_U$ будет возникать, как модуль коинвариантов $\overline{R}_{\mathcal{M}_n}=R/[R,R\mathcal{M}_n]$ модуля соотношений по действию $\mathcal{M}_n.$ Если рассматриваемое копредставление дискретно, то $M_U$ появятся, как $\overline{R}_{\mathcal{M}_n}=(R/[R,R\mathcal{M}_n])^{\wedge}_p.$ Использование $p$-адического пополнения в дискретном случае позволяет исключить в $M_U$ кручение взаимнопростое с $p$, а тогда квазирациональность влечет, что $\overline{R}_{\mathcal{M}_n}=(R/[R,R\mathcal{M}_n])^{\wedge}_p$ являются $\mathbb{Z}_p$- проективными $\mathbb{Z}_p[G/\mathcal{M}_n]$- модулями конечной $\mathbb{Z}_p$- размерности.
Если соответствующие проективные системы $\mathcal{T}_n, \mathcal{T}_n^p$ состоят, соответственно, из обобщенных пермутационных, пермутационных модулей, то мы будем говорить, что $\overline{R}=R/[R,R]$ - это обобщенный пермутационный  про-$\mathbb{Z}_pG^{\wedge}_p$-модуль, а $\overline{R}/p\overline{R}=R/R^p[R,R]$ - пермутационный про-$\mathbb{F}_pG^{\wedge}_p$-модуль.
Свойства пермутационности и обобщенной пермутационности модуля соотношений про-$p$ копредставления конечного типа (соответствующих проективных систем фактор- модулей коинвариантов в случае дискретных копредставлений) не зависят от выбора базиса окрестностей единицы для про-$p$-топологии свободной (про-$p$) группы, состоящего из нормальных делителей индекса, равного степени $p$, фильтрации Цассенхауза - лишь удобный выбор.

Обобщенный пермутационный модуль нельзя определить, как просто $\mathbb{F}_p$-пермутационный, поскольку из пермутационности модуля $R/R^p[R,R\mathcal{M}_n],$ вообще говоря, не следует отсутствие кручения в факторах коинвариантов $R/[R,R\mathcal{M}_n],$ требуемое в определении обобщенной пермутационности. Более того, конечная $p$-группа $G$ может действовать на абелевой группе $M=\mathbb{Z}/p^k\mathbb{Z}$ при $k\geq2$ нетривиально, поскольку например при $p>2$ имеем $\mid Aut(\mathbb{Z}/p^k \mathbb{Z}) \mid=(p-1)p^{k-1}$. Тогда $M$ - главный (то есть порожденный, как модуль, одним элементом) неразложимый модуль, поскольку он изоморфен циклической группе, но он не только не $\mathbb{Z}_p$-пермутационный модуль, но и не может иметь пермутационной системы образующих, как абелевой группы, в то время как $M/pM=\mathbb{Z} / p \mathbb{Z} $ уже является тривиальным пермутационным $\mathbb{F}_pG$-модулем.

Изощренность ситуации состоит в том, что для модулей из проективной системы $\mathcal{T}_n$ введенное Определение \ref{d1} является избыточным. На самом деле, при  $p>2$ обобщенная пермутационность совпадает с обычной  $\mathbb{Z}_p$-пермутационностью, а при $p=2$ возникает лишь действие на $\mathbb{Z}_2$ $"$изменением ориентации$"$, то есть с помощью автоморфизма аддитивной группы целых 2-адических чисел, переводящего $\pm 1$ в $\mp 1$. Действительно, в $\mathbb{Z}_p$ при $p>2$ нет корней $p$-й степени из единицы, поскольку отображение $log(1+z):1+p\mathbb{Z}_p=U^{(1)}\rightarrow \mathfrak{p}^1=p\mathbb{Z}_p$ устанавливает гомеоморфизм между мультипликативной группой главных единиц и аддитивной группой максимального идеала кольца $\mathbb{Z}_p$ \cite[Предл. 5.5]{Neu}. Но поскольку в $p\mathbb{Z}_p$ нет кручения, то и нет корней из единицы в мультипликативной группе кольца $\mathbb{Z}_p$. При $p=2$ имеется единственный - квадратный корень из единицы (равный $-1$). На самом деле, $-1$ является квадратным корнем из единицы в $\mathbb{Z}_2$. Остается проверить, что других корней из единицы в $\mathbb{Z}_2$ нет.  В этом случае имеем изоморфизм $1+2\mathbb{Z}_2\cong(1+4\mathbb{Z}_2)\times \{\pm1\}$ (логарифм сходится на $1+4\mathbb{Z}_2$ и задает гомеоморфизм $log(1+z):1+4\mathbb{Z}_2\rightarrow 4\mathbb{Z}_2$), а поэтому других корней из единицы нет.

Следующая Лемма является ключевым аргументом в доказательстве основного результата данной части работы и представляет из себя глубокий результат \cite[Теорема 3 и Следствие]{Weiss}.
\begin{lemma}\label{l1}
Пусть $M$ некоторый свободный $\mathbb{Z}_p[\zeta]$-модуль и задано представление конечной $p$-группы $G$ в $M$. Предположим, что $\overline{M}=M/\pi M,$ где $\pi= 1-\zeta,$ является $\mathbb{F}_p$- пермутационным $G$-модулем, тогда $M$ является обобщенным пермутационным $\mathbb{Z}_p[\zeta]G$- модулем. При этом из разложения $\overline{M}=\oplus_{i\in I}\mathbb{F}_p\uparrow^G_{H_i}$ следует разложение $M=\oplus_{i\in I}\mathbb{Z}_p[\zeta]\uparrow^G_{H_i}$.
\end{lemma}

\begin{theorem} \label{t0}
Пусть \eqref{eq1} - $QR$-(про-$p$) копредставление, тогда следующие свойства эквивалентны:

а) $\mathcal{T}_n$ - обобщенный пермутационный про-$G^{\wedge}_p$-модуль;

б) $\mathcal{T}_n^p$ -пермутационный про-$G^{\wedge}_p$-модуль.
\end{theorem}

\begin{proof}

а)$\Rightarrow$ б)
Поскольку $\mathbb{Z}_p[\zeta]$ вполне разветвлено над $\mathbb{Z}_p$.

б)$\Rightarrow$ а)
\emph{Доказательство проведем для про-$p$-копредставлений}, в дискретном случае рассуждения сохраняются, поскольку требуется лишь сравние $p$-адических пополнений фактормодулей коинвариантов с их $mod(p)$- факторами (откуда и возникают про--объекты в формулировке теоремы).
Так как $\overline{R}/p\overline{R}$ - пермутационный модуль, то $\overline{R}/p\overline{R}\cong \oplus_{i\in I} \mathbb{F}_p\uparrow^G_{H_i}.$ В дискретном случае мы точно знаем лишь, что $|I|=dim_{\mathbb{F}_p}R/R^p[R,F]$, в про-$p$- случае из \cite[Теоремы 3.2]{Mel1} выукает, что $|I|=dim_{\mathbb{F}_p}H^2(G,\mathbb{F}_p),$ а $H_i$ - это конечные циклические подгруппы в $G$.
Если мы докажем, что для каждого $n$ из пермутационности модуля $(\overline{R}/p\overline{R})_{\mathcal{M}_n}=R/R^p[R,R\mathcal{M}_n]\cong\oplus_{i\in I} \mathbb{F}_p\uparrow^{G/\mathcal{M}_n}_{H_i\mathcal{M}_n/\mathcal{M}_n}$
 следует разложение $R/[R,R\mathcal{M}_n]\cong\oplus_{i\in I} \mathbb{Z}_p\uparrow^{G/\mathcal{M}_n}_{H_i\mathcal{M}_n/\mathcal{M}_n}$,
то утверждение будет доказано, но это содержание Леммы \ref{l1}. Таким образом получаем
$R/[R,R\mathcal{M}_n] \cong \oplus_{i\in I} \mathbb{Z}_p\uparrow^{G/\mathcal{M}_n}_{H_i\mathcal{M}_n/\mathcal{M}_n},$
где в правой части стоит прямая сумма главных обобщенных пермутационных модулей (Определение \ref{d01}).
\end{proof}

Напомним, $\mathfrak{p}$-модулярной системой \cite[1.9]{Be} называется тройка $(K_{\mathfrak{p}},\mathcal{O}_{\mathfrak{p}},k),$ где $\mathcal{O}_{\mathfrak{p}}$ -полное кольцо дискретного нормирования ранга 1 ($\mathcal{O}_{\mathfrak{p}}$ - область главных идеалов с единственным ненулевым простым идеалом) с максимальным идеалом $\mathfrak{p}=(\pi)$ и полем вычетов $k$ характеристики $p>0$, а $K_{\mathfrak{p}}$ - его поле частных характеристики 0. Использование про-$p$-пополнения при построении $\mathcal{T}_n=((\frac{R} {[R,R\mathcal{M}_n]} )^{\wedge}_p, \phi^{n+1}_{n}),\mathcal{T}_n^p=(\frac{R} {R^p[R,R\mathcal{M}_n]}, \widetilde{\phi}^{n+1}_{n})$ в дискретном случае помимо исключения кручения взаимнопростого с $p$ позволяет воспользоваться преимуществами $\mathfrak{p}$-модулярной системы $(\mathbb{Q}_p,\mathbb{Z}_p,\mathbb{F}_p)$
в качестве коэффициентов представлений. На самом деле, по теореме Грина $"$о неразложимости$"$ \cite[3.13.3]{Be} модули $\mathbb{F}_p \uparrow^{G/\mathcal{M}_n}_{H_i\mathcal{M}_n/\mathcal{M}_n}$ являются абсолютно неразложимыми (неразложимы над любым расширением поля $\mathbb{F}_p$).
Более того для $\mathbb{Z}_p$, как  кольца коэффициентов представления, имеет место теорема Крулля-Шмидта \cite[1.4]{Be}, что вместе с теоремой Грина \cite[3.10.2]{Be} влечет жесткость разложения $R/[R,R\mathcal{M}_n] \cong \oplus_{i\in I} \mathbb{Z}_p\uparrow^{G/\mathcal{M}_n}_{H_i\mathcal{M}_n/\mathcal{M}_n}$ с точностью до порядка классов изоморфизмов факторов редукции по модулю $p$.

\begin{sled} \label{H1}
Класс квазирациональных про-$p$-копредставлений обладает следующими свойствами:

1) асферические про-$p$-копредставления и копредставления про- $p$-групп с одним соотношением являются квазирациональными;

2) модуль соотношений $\overline{R}=R/[R,R]$ квазирационального копредставления является обобщенно пермутационным тогда и только тогда, когда $\mathbb{F}_p$- пермутационным является его $mod(p)$ фактор $R/R^p[R,R]$.

\end{sled}

\begin{proof}

Нами уже отмечалось, что квазирациональными являются про-$p$- копредставления про-$p$-групп с одним соотношением, копредставления и подкопредставления $CA$- асферических копредставлений, а также асферические про-$p$-копредставления. Таким образом квазарациональность удовлетворяет условию 1) Следствия \ref{H1}. Из Теоремы \ref{t0} следует обобщенная эквивалентность $\mathbb{Z}_p$ и  $\mathbb{F}_p$- пермутационности, а поэтому $QR$-(про-$p$) копредставления  удовлетворяют условию 2) Следствия \ref{H1}.
\end{proof}

Мы оставляем читателю самостоятельно сформулировать и доказать аналог Следствия \ref{H1} для дискретных $QR$-копредставлений в терминах про-объектов. В этом случае условие 1) Следствия привлекательно конструировать используя $CA$-асферические копредставления и их подкопредставления.

В самом начале работы \cite[стр.71]{Mel1} О.В. Мельников отмечает, что существует теория $\mathbb{Z}_p$-пермутационных модулей соотношений полностью параллельная теории асферических про-$p$-копредставлений (то есть когда модули соотношений являются $\mathbb{F}_p$- пермутационными). Также в личных беседах с автором еще в 1997 году он неоднократно высказывал неудовлетворенность отсутствием способов выявления асферичности копредставления в ряде ключевых ситуаций (например не получается выяснить всякая ли про-$p$-группа с одним соотношением является асферической).
Мельников предполагал существование некоторого класса про-$p$-копредставлений (с содержательным аппаратом изучения и возможностью достаточно просто проверить принадлежность к этому классу) со следующими свойствами:

1) данный гипотетический класс содержит одновременно асферические про-$p$- копредставления и копредставления про-$p$-групп с одним соотношением;

2) в рамках данного гипотетического класса копредставлений возможна унификация концепций $\mathbb{Z}_p$ и $\mathbb{F}_p$-пермутационности.

В работах \cite{Mikh2014,Mikh2015,Mikh2016} показана возможность исследования квазирациональных копредставлений с помощью перехода к $p$-адическим Мальцевским пополнениям. Таким образом, Следствие \ref{H1} подтверждает правильность предположения, высказанного О.В. Мельниковым.


\end{document}